\documentclass[lettersize,journal]{IEEEtran}
\usepackage{amsmath,amsfonts}
\usepackage{algorithmic}
\usepackage{array}
\usepackage[caption=false,font=normalsize,labelfont=sf,textfont=sf]{subfig}
\usepackage{textcomp}
\usepackage{stfloats}
\usepackage{url}
\usepackage{verbatim}
\usepackage{graphicx}
\usepackage{caption}
\def\BibTeX{{\rm B\kern-.05em{\sc i\kern-.025em b}\kern-.08em
    T\kern-.1667em\lower.7ex\hbox{E}\kern-.125emX}}
\usepackage{balance}
\usepackage{tikz}
\usepackage{xcolor}
\usetikzlibrary{patterns}
\usepackage{algorithm}
\usetikzlibrary{arrows.meta,positioning}
\usepackage{thmtools}
\usepackage{thm-restate}
\usepackage{amssymb, amsthm}
\usepackage{microtype}
\usetikzlibrary{decorations.pathreplacing}
\usepackage{subfig}
\usepackage{url}

\usepackage{pgfplots}
\pgfplotsset{compat=1.17}
\usepgfplotslibrary{patchplots}
\usepgfplotslibrary{statistics}
\usepgfplotslibrary{external} 
\usepgfplotslibrary{fillbetween}
\usepackage{subfig}
\usetikzlibrary{external}
\usepackage{pgfplotstable}
\pgfplotsset{select coords between index/.style 2 args={
    x filter/.code={
        \ifnum\coordindex<#1\fi
        \ifnum\coordindex>#2\fi
    }
}}

\definecolor{colorpyC0}{RGB}{31, 119, 180}
\definecolor{colorpyC1}{RGB}{255, 127, 14}
\definecolor{colorpyC2}{RGB}{44, 160, 44} 
\definecolor{colorpyC3}{RGB}{214, 39, 40} 
\definecolor{colorpyC4}{RGB}{148, 103, 189}
\definecolor{colorpyC5}{RGB}{140, 86, 75} 
\definecolor{colorpyC6}{RGB}{227, 119, 194}
\definecolor{colorpyC7}{RGB}{127, 127, 127}
\definecolor{colorpyC8}{RGB}{188, 189, 34} 
\definecolor{colorpyC9}{RGB}{23, 190, 207}

\newcommand{\focs}{\textsc{Focs}}

\newcommand{\cotwo}{\textsc{CO}_2}
\newcommand{\costfunction}{\cotwo}

\def\BibTeX{{\rm B\kern-.05em{\sc i\kern-.025em b}\kern-.08em
    T\kern-.1667em\lower.7ex\hbox{E}\kern-.125emX}}

\begin{document}
\title{Carbon, Cost and Capacity: Multi-objective Charging of Electric Buses}
\author{Leoni Winschermann, Leander C. van der Bijl, Marco E.T. Gerards and Johann Hurink
\thanks{Submitted on 08.04.2025.
This research is conducted within the \textit{MegaMind} program subsidized (partly) by the Dutch Research Council (NWO) through the Perspectief program under number P19-25 and the \textit{SmoothEMS met GridShield} project subsidized by the Dutch ministries of EZK and BZK (MOOI32005). The first and second author share main authorship and flipped a coin to determine order of authors.}
\thanks{All authors are with the Department of Electrical Engineering, Mathematics and Computer Science,
University of Twente, 7522NB Enschede, The Netherlands (corresponding authors: l.winschermann@utwente.nl, l.c.vanderbijl@utwente.nl) }
}

\maketitle

\begin{abstract}
The public transport sector is in the process of decarbonizing by electrifying its bus fleets. This results in challenges if the high electricity demand resulting from battery charging demand is confronted with limited grid capacity and high synchronicity at bus charging sites. 
In this paper, we explore multi-objective scheduling for bus charging sites to minimize the emissions associated with charging processes and to aid the operation of the electricity grid by mitigating peak consumption. 
In particular, we discuss and validate optimization approaches for those objectives, as well as their weighted combination, based on data from a real-life bus charging site in the Netherlands. 
The simulation results show that compared to uncontrolled charging, power peaks can be reduced by up to 57\%, while time-of-use emissions associated with the charging of electric buses are also reduced significantly. 
Furthermore, by using a synthetic baseload, we illustrate the flexibility potential offered by bus charging sites, and advocate that such sites should share a grid connection with other high-load assets.
\end{abstract}


\section{Introduction}
\label{bus:sec:intro}
\IEEEPARstart{T}{he} Dutch electricity grid currently struggles with severe congestion \cite{hofstedeStroomnettenZijnOvervol2023}. This, among others, prevents industrial parties from obtaining a (larger) grid connection, while the pressure to electrify their processes increases due to the energy transition and the ambition to become climate neutral. Furthermore, companies can even lose their capacity during peak hours due to new legislation \cite{ToezichthouderACMWil2024}. These restrictions may enforce sudden changes in operation with possibly high costs as a result.

One of the industries that must electrify their processes is the public transport sector. As a consequence, diesel buses are being replaced by electric buses. Their batteries must be charged, which usually takes place at a charging site that charges multiple buses at once. The magnitude of the resulting synchronous energy charging demands, in combination with the high charging rates bus batteries support, results in high peak powers if the charging occurs uncontrolled. Furthermore, the main portion of the charging occurs at hours just after the individual buses seize operating. Given limited-size bus fleets, this most likely will be during the night, which especially in summer is the time of day with high emissions associated to the electricity drawn from the grid. Therefore, a control method for the charging of the buses should take both grid stress and $\cotwo$ emissions into account. In this context, it is important to also investigate if the bus company can coordinate its electricity usage with other companies with buildings closeby to ensure stable operation and make more efficient use of the available grid capacity and (locally generated) renewable energy.

This article considers the real-life scenario of a bus charging site owned by a public transport company in the Netherlands. This company is scaling up the number of electric buses to comply with regulations and because of their ambition to reduce their $\cotwo$ emissions to zero as soon as possible. In this article we present a bus charging algorithm that addresses the aforementioned problems and numerically evaluate the performance of this algorithm.

The topic of charging electric buses has only sparingly been treated in literature. The overview in \cite{perumalElectricBusPlanning2022} surveys 43 articles related to electric bus charging and scheduling and identifies three stages with respect to this topic, namely:
\begin{itemize}
    \item Strategic planning (time span: years). This includes the planning of the bus lines and charging infrastructure.
    \item Tactical planning (time span: one year to days). This includes vehicle scheduling and charging scheduling.
    \item Operational planning (time span: intraday to real-time).
\end{itemize}

The research presented in this paper focuses on charging schedules, which can be created based on the given driving schedule of the buses. Therefore, the work contributes to the tactical planning stage as defined above. Out of the 43 papers reviewed in~\cite{perumalElectricBusPlanning2022}, only 11 of them optimize the charging schedule, most with the sole objective to minimize electricity (or battery) costs. 

This topic is already well researched for EVs in general (see, e.g.~\cite{2023KleinEVchargeSchedwithFlexibleServiceOperations,zhengOnlineDistributedMPCBased2019,shindeOptimalElectricVehicle2016,alonsoOptimalChargingScheduling2014}). These papers present charging algorithms for minimizing time-of-use costs or for flattening of the total profile. In these cases, uncertainty is an important aspect as it is often unknown beforehand when cars will arrive and leave. Even when these times are known, still a proper incentive has to be given to users to actually follow the solutions provided by such algorithms (see e.g. \cite{adetunjiTwoTailedPricingScheme2024}).

In the bus charging case, the bus schedule is already known far in advance and therefore the behaviour is more predictable. However, it is still possible that deviations from the schedule occur due to delays, implying that there is some uncertainty. One can ignore these uncertainties if the dwell time of a bus is at the charging site is large enough (for example 18 hours) with respect to the arrival time uncertainty (usually within half an hour), as is the case for the considered situation. As mentioned above, the most commonly considered objective within the deterministic bus charging context is cost minimization, see e.g.~\cite{zhuang_stochastic_2021,liCooperativeOptimizationBus2023,panOptimizationElectricBus2020,basmaOptimizationBatteryElectric2020,chengweizhangChargingScheduleOptimization2019,leouOptimalChargingSchedule2017,zhouCollaborativeOptimizationVehicle2020,abdelwahedEvaluatingOptimizingOpportunity2020,heOptimalChargingScheduling2020a,bagherinezhadSpatioTemporalElectricBus2020}. In these studies the goal is to minimize time-of-use costs while taking grid limits into account. Note that these time-of-use cost minimization algorithms can also be adopted to minimize $\cotwo$ emissions as one can simply use the emissions of the electricity mix as a `price' input. However, it is in general not possible for these algorithms to optimize with respect to a weighted combination of costs and flattening of the energy profile. On the other hand, Jahic \textit{et al.} \cite{jahicChargingScheduleLoad2019} consider flattening of the aggregated power profile of a bus depot by minimizing the sum of squares of the power profile. However, it is not possible to simultaneously minimize $\cotwo$ emissions using their approach, and they consider a non-preemptive problem formulation, i.e., once a bus started charging, it keeps charging till full. On top of that, both algorithms they propose are heuristics. 

Based on the above, in this paper we aim to combine these two objectives. Given a bus schedule for a day, this paper contributes:
\begin{itemize}
    \item An open-source algorithm that solves the bus charging scheduling problem for a weighted combination of time-of-use costs (in this paper $\cotwo$ emissions are used) and profile flattening.
    \item A numerical evaluation of the effect of the chosen weights on either objective value.
    \item A validation of the approach with a real-life bus schedule from a charging site in the Netherlands.
    \item An analysis of the flexibility potential of bus charging hubs to flatten out large baseloads.
\end{itemize}

The remainder of the paper is structured as follows.
Section~\ref{sec:busproblemDescr} formalizes the considered problem, discusses three objective functions and presents solution approaches for each of them. Next, Section~\ref{sec:method:simstudy} describes
a simulation study, the used data sets,
and the considered evaluation metrics
. The resulting power profiles for the use case are discussed in Section~\ref{bus:sec:results} together with the observed trade-offs. Specifically, Section~\ref{sec:results:baseloadflex} highlights the flexibility potential of bus charging sites.
Finally, Section~\ref{bus:sec:conclusion} summarizes the findings of the paper.

\section{Methods}\label{bus:sec:methods}
In this section, we state the problem definition (Section~\ref{sec:busproblemDescr}) and discuss the three objective functions used in this paper: $\cotwo$ minimization
, profile flattening expressed as sum of squares minimization 
and a weighted combination of the two
. For each objective, we give its formal definition and describe a corresponding optimization approach. Next, Section~\ref{sec:method:simstudy} describes the setup for the simulation study, the real-world case data and the considered evaluation metrics.

    \subsection{Problem description} \label{sec:busproblemDescr}
        Compared to other electric vehicle scheduling problems, bus charging is strongly dependent on a timetable specifying trips per line and bus. Therefore, the bus charging problem holds quite a high planability, which allows to consider the day-ahead optimization problem as a deterministic problem. In this section, we introduce the notation, model constraints and various considered objective functions.

         We define the following scheduling problem for charging the buses. 
        Given are a planning horizon discretized into unit-sized time intervals $M_i$, $i\in \mathcal{M} := \{1,\dots,m\}$ and a set of to be scheduled charging sessions, denoted by $\mathcal{J} := \{1,\dots,n\}$. Each job $j\in\mathcal{J}$ has an associated arrival time $a_j$, a departure time $d_j$, an energy requirement $e_j$ and a maximum charging rate $\ell_j$. The latter is the maximum amount of energy that can be charged to job $j$ within a unit-sized time interval. We assume arrival and departure times to intersect with the start and end points of intervals in $\mathcal{M}$.
        A job $j$ is said to be available in an interval $M_i$ if its arrival time is before the start of the interval, and its departure is no earlier than the end of the interval. 
        Let $J(i)$ be the set of buses available for charging in interval $M_i$, $i\in \mathcal{M}$. 
        Similarly, let $J^{-1}(j)$ be the set of indices $i$ of intervals $M_i$ where job $j \in \mathcal{J}$ is available.

        To derive a charging schedule, we introduce decision variables $e_{i,j}$, which express the amount of energy charged to job $j$ during interval $M_i$. We only introduce $e_{i,j}$ for jobs $j$ which are available during $M_i$.
        The constraints can be formalized as follows:
        \begin{subequations}\label{eq:busMIP}
        \begin{align}
            \sum_{i\in J^{-1}(j)} e_{i,j} &\geq e_{j}  &\forall j\in\mathcal{J}\hspace{2.5pt} \label{eq:busMIPnoENS}\\
            e_{i,j} &\geq 0& \forall j \in \mathcal{J}, i\in J^{-1}(j)\hspace{2.5pt}  \label{eq:busMIPnonnegLoad}\\
            e_{i,j} &\leq \ell_j & \forall j \in \mathcal{J}, i\in J^{-1}(j). \label{eq:busMIPspeedlimit}
        \end{align}
        \end{subequations}
        Let $s$ be the aggregated power profile of a solution. Given that we consider unit-sized intervals, this can be normalized to energy, leading to $s(i) = \sum_{j\in J(i)} e_{i,j}$.

        As mentioned before, two of the most frequently discussed objectives in energy management are price steering and a minimization of the sum of squares of the (discretized) power profile, or intuitively the flattening of the power profile. We consider them one by one, and then their combination in a weighted objective function. 
        
        \subsubsection{$\cotwo$ minimization}\label{sec:method:obj:cotwo}
            First, for minimization of $\cotwo$ emissions associated with bus charging, we focus on the emission factors of the electricity mix. Let $\costfunction(i)$ be the emission factor in kg$\cotwo{eq}/$kWh associated with the electricity mix during interval $M_i$. This leads to the following notation for the objective function:
            \begin{align}
                C(s) = \sum_{i\in\mathcal{M}}s(i) \costfunction(i)\label{eq:costobjectiveDiscretized}
            \end{align}

\begin{figure}[]
                \centering
                \begin{tikzpicture}[scale=1, every node/.style={scale=1}, 
                  mycircle/.style={
                     circle,
                     draw=black,
                     fill=gray,
                     fill opacity = 0.2,
                     text opacity=1,
                     inner sep=0pt,
                     minimum size=20pt,
                     font=\small},
                  myarrow/.style={-Stealth},
                  node distance=0.6cm and 1.2cm
                  ]      
                  \node[mycircle] (j1) at (1,1) {$1$};
                  \node[mycircle] (j2) at (1,-0.25) {$2$};
                  \node[rotate=90] (jdots) at (1,-1) {\dots};
                  \node[mycircle] (jn) at (1,-2) {$n$};
            
                  \node[mycircle] (s) at (-1,-0.5) {$v_0$} ;
            
                  \node[mycircle,] (i1) at (3.25,1){$M_1$};
                  \node[mycircle] (i2) at (3.25, -0.25){$M_2$};
                  \node[rotate=90] (idots) at (3.25,-1) {\dots};
                  \node[mycircle] (im) at (3.25, -2) {$M_m$};
                  
                  \node[mycircle] (t) at (5.25,-0.5) {$v_t$};
            
                \foreach \i/\j/\txt/\p in {
                  s/j1/$e_1$/above,      
                  s/j2/$e_2$/above,
                  j1/i1/$\ell_1$/above,      
                  i1/t/$g_1$/above   
                  }
                   \draw [myarrow] (\i) -- node[font=\small, sloped, above] {\txt} (\j);
            
                \foreach \i/\j/\txt/\p in {
                  s/jn/$e_n$/below,
                  j2/i2/$\ell_2$/below, 
                  jn/im/$\ell_n$/below,
                  i2/t/$g_2$/below,
                  im/t/$g_m$/below
                  }
                   \draw [myarrow] (\i) -- node[font=\small, sloped, below] {\txt} (\j);
                   \draw [myarrow] (j1) -- node[font=\small, sloped, above=-2pt, pos = 0.35] {$\ell_1$} (i2);
                   \draw [myarrow] (j2) -- node[font=\small, sloped, below=-1.7pt , pos = 0.35] {$\ell_2$} (i1);
                \end{tikzpicture}
                \caption{Schematic of flow network structure for the $\cotwo$\ minimization for charging electric buses.}
                \label{fig:mincostflowNetwork}
            \end{figure}
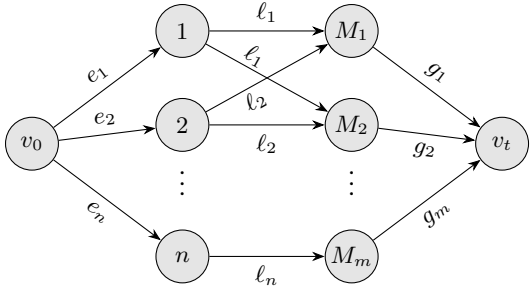

            For this objective, Problem~\eqref{eq:busMIP} is equivalent to a minimum cost flow problem (see Fig.~\ref{fig:mincostflowNetwork} for a suitable network structure to model the problem).
            Next to sink and source nodes $v_0$ and $v_t$, there is a node layer corresponding to the job set $\mathcal{J}$, and a node layer corresponding to the unit-sized intervals $\mathcal{M}$ partitioning the time horizon. The annotations at the edges are the respective edge capacities. 
            In particular, the edge layer between job and interval nodes models the schedule. Given a suitable flow, the flow through the edge between job node $j$ and interval node $M_i$ corresponds to decision variable $e_{i,j}$, which expresses the energy scheduled for $j$ during $M_i$. To ensure compliance with \eqref{eq:busMIPspeedlimit}, a job node $j$ is connected to interval node $M_i$ if and only if $j\in J(i)$. In that case, the maximum charging rate $\ell_j$ of the job limits the flow capacity through that edge. 

            Moreover, all interval nodes are connected to sink node $v_t$. By flow conservation, the flow outgoing from an interval node $M_i$ equals the sum of the incoming flows, which in turn model the decision variables $e_{i,j}$. Optionally, one may introduce global capacity constraints $g_i$ for intervals $M_i$ on those edges. To ease the presentation, we consider $g_i = \infty$, but the presented concepts can be extended to finite bounds. 

            Finally, there is an edge between source node $v_0$ and each job $j$ with capacity $e_j$. Intuitively, we are interested in flows that saturate all edges out-going from $v_0$ so that \eqref{eq:busMIPnoENS} is tight for every job $j$. By limiting ourselves to positive flows, we satisfy constraint \eqref{eq:busMIPnonnegLoad}. 

            To derive a minimum cost flow, we associate the edge between interval node $M_i$ and sink $v_t$ with a cost corresponding to $\costfunction(i)$ for all $i\in \mathcal{M}$. All other edge costs are set to 0. Solving a minimum cost flow with flow value $\sum_j e_j$ gives a schedule minimizing $\cotwo$ emissions \eqref{eq:costobjectiveDiscretized} subject to constraints \eqref{eq:busMIP}.

        \subsubsection{Sum of squares minimization} \label{sec:method:obj:flat}
            As a second objective, we consider the minimization of the squares of the (discretized) power profile, to express the flattening of the aggregated power profile. This objective corresponds to the square of the 2-norm of the aggregated power profile, or formally
            \begin{align}
                F(s) = \sum_{i\in \mathcal{M}} s(i)^2. \label{eq:flatobjectiveDiscretized}
            \end{align}
             We can use the open source solver \focs ~\cite{2023WinschermannFOCSArxiv} to solve Problem~\eqref{eq:busMIP} with objective \eqref{eq:flatobjectiveDiscretized}. In fact, \focs\ minimizes a wide range of objective functions under constraints \eqref{eq:busMIP}. Given a fixed time discretization, it finds the optimal solution in time polynomial in the number of buses $n$.

        \subsubsection{Weighted combination of objectives} \label{sec:method:obj:weighted}
            As third objective, we introduce the weighted combination of both objectives
            \begin{align}
                W(s,w_c,w_f) = w_cC(s) + w_fF(s), \label{eq:weightedobjectiveDiscretized}
            \end{align}
            where $w_c$ and $w_f$ are strictly positive weights for the $\cotwo\ $ and flatness objectives respectively.

            The reason for using a weighted combination of both objectives is that they both consider important but essentially different aspects and using one of either could lead to improving one objective while significantly worsening the other. In Section \ref{bus:sec:results} we show that this happens in the practical case considered in this paper. It therefore makes sense to use a combination of the two to ensure both objectives are considered.

            As the objectives use two different units (kW and $\cotwo{\text{eq}}$), $w_f = w_c$ does not mean that both objectives will be taken evenly into account. The only way to see what a sensible balance of both weights is, is to first solve the problem with unit weights and to adapt then based on the achieved energy profile. A figure with results for different ratios is given in Section \ref{bus:sec:results}.
            
            We may solve the problem with the weighted objective function using \focs. However, as stated before, \focs\ only minimizes a certain class of objective functions. Therefore, we have to transform weighted the objective $W(s,w_c,w_f)$ to a form from which we derive auxiliary jobs to add to the input instance such that \focs\ minimizes \eqref{eq:weightedobjectiveDiscretized} under constraints \eqref{eq:busMIP}.
            
            First, note that the minimum of \eqref{eq:weightedobjectiveDiscretized} can be written as 
            \begin{align*}
                \min \ \sum_{i\in\mathcal{M}} s(i)^2 + \frac{w_c}{w_f} s(i) \costfunction(i).
            \end{align*}
            Our aim is to rewrite for each $i\in \mathcal{M}$ the corresponding term in the form 
            \begin{align*}
                (\alpha s(i) + \beta)^2 - \gamma .
            \end{align*}
            This yields 
            \begin{align*}
                \alpha & = 1,\\
                \beta & = \frac{w_c}{2 w_f} \costfunction(i),\\
                \gamma & = \beta^2 = \frac{w_c^2}{4w_f^2} \costfunction(i)^2.
            \end{align*}
            Note that $\gamma$ is a constant which is independent of $s$. Therefore, we can consider the equivalent objective function 
            \begin{align}
                \min \ \sum_{i\in \mathcal{M}} \left(s(i) + \frac{w_c}{2w_f}\costfunction(i)\right)^2. \label{eq:objectiveMinCostTransformed}
            \end{align}
            
            To model this objective in \focs\ we add $\frac{w_c}{2w_f}\cotwo (t)$ as a baseload of auxiliary jobs
            . 
            Note that instead of $\costfunction(i)$ we may also consider a weighted combination of financial and environmental costs
            . 
        
    \subsection{Simulation study}\label{sec:method:simstudy}
        In this paper we focus on simulating the operation of a bus charging site over the course of one week when optimizing for the various objective functions described above.\footnote{The code used in the simulations is available under \url{https://github.com/lwinschermann/FlowbasedOfflineChargingScheduler}} All experiments are based on real-world data, and schedules are computed for 15-minute granularity. In the following, we describe the used data in more detail.
        
        \subsubsection{Data set} \label{sec:busData}
            A bus line schedule is used from a real-life bus charging site in the Netherlands
            . The data set includes for each day of the week all lines of buses that charge at that charging site. Per line, the following is given:

            \begin{itemize}
                \item The start and end times of the line.
                \item The type of bus needed for the line.
                \item The state of charge of the bus after driving the line.
            \end{itemize}
            
            There are two types of buses in this data set: one with a 122~kWh battery and one with a 273~kWh battery. The number of lines for a day depends on the day of the week. Specifically, for Monday until Thursday the same schedule is used, whereas Friday, Saturday and Sunday each have their own schedule.

            Table~\ref{tab:numberbusesperday} shows the number of lines per day. Note that there is an equal number of lines for Monday up to Friday but Saturday and Sunday have less lines.

            \begin{table}[h]
            \centering
            \caption{Statistics per day of the week.}
            \label{tab:numberbusesperday}
            \begin{tabular}{l|c c c}
            \hline
            & Monday - Friday & Saturday & Sunday \\
            \hline
            \# lines per day & 33 & 22 & 23\\
            \# matched buses per day & 33 & 22 & 21\\
            \hline
            \end{tabular}
            \end{table}

            On the charging site, there are sufficient $30$~kW chargers available to charge every bus. We assume that every bus has to be charged to 100\% state of charge before starting a line.
            
            To determine by what time a bus must be fully charged again to drive a line the next day, we create an assignment of the buses arriving each evening to the lines of the next day. This is done by creating a bipartite graph where the nodes $B$ on one side represent buses and the nodes $L$ on the other side represent lines for the next day. An edge $(b,l)$, $b\in B$, $l \in L$, exists if there is sufficient time for bus $b$ to be fully charged before the starting of line $l$ based on the required capacity and the charging limit of $30$~kW. A maximal matching is then calculated resulting in the bus-line allocation that will be used the next day. Buses that cannot be matched receive the end of the next day as deadline for their charging. This matching simplifies the problem and aims to limit the total number of buses needed to drive all lines.

            The number of buses that can be matched to a line in the following day can be seen in the second row of Table~\ref{tab:numberbusesperday}. When comparing this to the first row (the number of buses needed on the respective days) it can be seen that on almost every day, all used buses can be matched to a line for the next day. Only on Sunday two buses cannot be matched and therefore receive the end of the next day as charging deadline. Note that this implies that two additional spare buses will be scheduled for the next day.

            The bus charging site shares a 950~kW grid connection with a small office building that has a solar rooftop installation. For this research we used available data describing the energy drawn from the grid to supply the office building
            per $15$~minute interval. We treat this load as static baseload for the charging schedule. 
            For the $\cotwo$ emission factors, datasets in 
            hourly time granularity are available under \cite{2024nednlData}. (This data is from an initiative by the Dutch transmission grid operators for gas and electricity, where they shared a dataset with two weeks of emission factors in 15 minute granularity).

        \subsubsection{Evaluation metrics} \label{bus:sec:evalmetrics}
        We evaluate the simulation outputs on the following three metrics. 
        Given a schedule $s$, the first two natural metrics for the given setting in this paper are the $\cotwo$ emissions $C(s)$ per unit of energy and the flatness $F(s)$ of the aggregated power profile as defined in Equations~\eqref{eq:costobjectiveDiscretized} and~\eqref{eq:flatobjectiveDiscretized} respectively.

        In the Netherlands, peak tariffs and the reduction of grid connections are of increasing interest for industrial parties. Furthermore, bus fleets are in the process of electrification, implying that the number of electric buses and hence the total energy demand is still increasing, having implications for the grid connection capacity needed to feasibly operate a bus charging site in the future. 
        Therefore, as a third metric we consider the peak power of a given schedule. This is defined as 
        \begin{align}
            P(s) = \max_{i\in \mathcal{M}} \frac{s(i)}{|M_i|}, \nonumber
        \end{align}
        where $|M_i|$ is the length of interval $M_i$. 
        
\section{Results} \label{bus:sec:results}
This section discusses the results of the simulation study. 
All results are based on baseload data corresonding to June $5-11$ of the year $2023$, as well as $\cotwo$ emission data from June 1 2023 onwards.
To this end, we present the power profiles of the bus charging hub for both an uncontrolled scenario and for the various objective functions discussed in Section~\ref{sec:busproblemDescr}. Furthermore, we consider the trade-off between flattening and $\cotwo$ objectives in the weighted scenario in more detail. 
Finally, Section~\ref{sec:results:baseloadflex} illustrates the potential flexibility that a bus charging site holds. 

\subsection{Bus charging schedules}

\begin{figure*}
    \centering
    \includegraphics[]{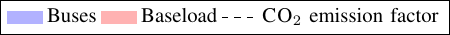}
    
    \subfloat[The Power Profiles for uncontrolled charging.]{\label{fig:bus:powerprofiles_noglobalcal_greedy}
    \includegraphics[]{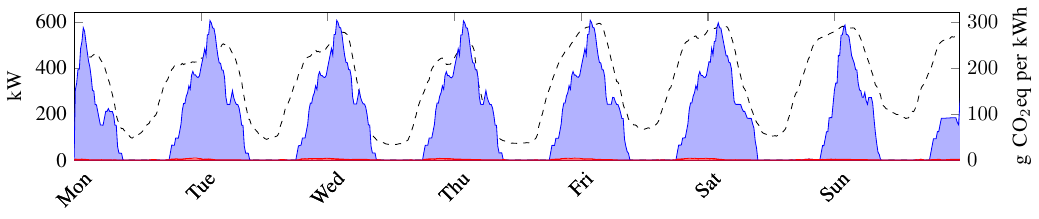}
    }
    \\
    \subfloat[The Power Profiles for CO$_2$ steering.]{\label{fig:bus:powerprofiles_noglobalcal_cotwo}
    \includegraphics[]{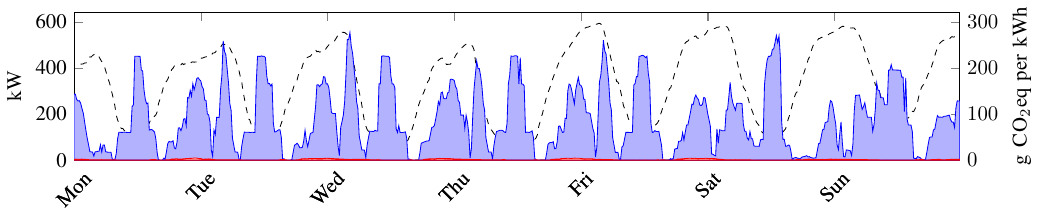}
    }
    \\
    \subfloat[The Power Profiles for $w_c = 1$ and $w_f = 2$.]{\label{fig:bus:powerprofiles_noglobalcal_weighted}
    \includegraphics[]{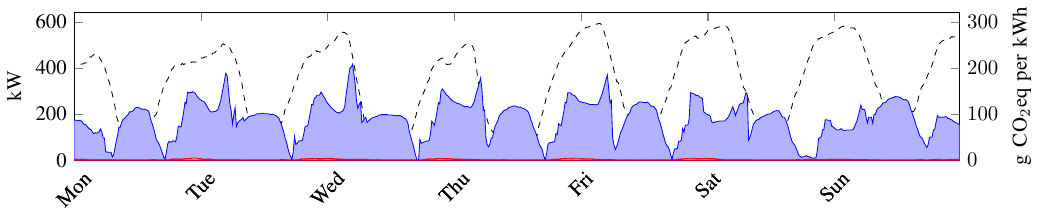}
    }
    \\
    \subfloat[The Power Profiles for flattening steering.]{\label{fig:bus:powerprofiles_noglobalcal_flat}
    \includegraphics[]{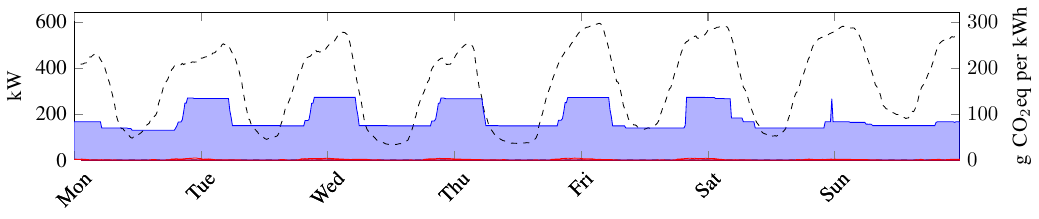}
    }
    \\
    \caption{The power profiles with data-based baseload.}
    \label{fig:bus:powerprofiles_noglobalcap}
\end{figure*}

To start with, Fig.~\ref{fig:bus:powerprofiles_noglobalcap} shows the aggregated power profiles of the simulated week at the bus charging site. The data is plotted in 15 minute intervals. The left-hand y-axis represents the aggregated power of the charging site and corresponds to the power curves associated with the baseload of the adjacent building (red plot) and the power used to charge buses (blue plot). In the figure, these two power demands are stacked on top of each other. We display a total of four scenarios: (\ref{fig:bus:powerprofiles_noglobalcal_greedy}) is the uncontrolled scenario where buses charge at their maximum charging power starting from their arrival until they are fully charged. Figures~\ref{fig:bus:powerprofiles_noglobalcal_cotwo} and \ref{fig:bus:powerprofiles_noglobalcal_flat} correspond to scheduling with objective functions $C(s)$ and $F(s)$ respectively, and Fig.~\ref{fig:bus:powerprofiles_noglobalcal_weighted} corresponds to the weighted objective with $w_c = 1$ and $w_f = 2$. These values are arbitrary and serve as an example of what such a balance looks like.
Lastly, the dashed plot depicts the $\cotwo$ emission factor associated with the respective time interval and the corresponding units are given on the right-hand y-axis.

As mentioned in Section~\ref{bus:sec:intro}, the bus charging site considered in this case study is currently electrifying their bus fleet. In particular, the amount of buses may double over the next years. While the current simulation results show that even in the uncontrolled scenario (see Fig.~\ref{fig:bus:powerprofiles_noglobalcal_greedy}) the power peaks do not exceed the grid-connection capacity constraint of 950~kW given the current number of electric buses, this will not apply once the entire fleet runs on electricity. In fact, the bus charging site already reports grid capacity violations on for example holidays when more shifts finish early in the evening and the synchronicity of charging is higher. Therefore, the power peak reductions observed in the controlled cases become even more meaningful, not just in terms of grid connection tariff savings, but in keeping the charging site operational.

Another observation is that the baseload (red plot) is marginal given the magnitude of the bus charging demand. This is mainly due to the relatively small size of the adjacent office building. Its highest power consumption reported in 2023 is 12.3~kW, and in the considered week in early June there was additional solar generation covering a large portion of the baseload. To investigate the potential flexibility of a bus charging hub like this, in Section~\ref{sec:results:baseloadflex} 
we investigate the flattening potency for larger (synthetic) baseloads.

Considering the various objective functions (see Table~\ref{tab:busprofilesObjectiveFunctions} for the objective values per scenario), we first analyze the base case scenario: uncontrolled charging.  
\begin{table}[h]
\centering
\caption{Evaluation metrics for bus charging schedules over a full week under various scenarios.}
\label{tab:busprofilesObjectiveFunctions}
\begin{tabular}{l|rrr}
\hline
  & $F(s)$ & $ C(s)$ & $ P(s)$ \\
       & $[10^{7}]$& $[10^{6}~\text{kg}\cotwo\text{eq}$]  & [kW]  \\
\hline
Uncontrolled  & 4.37 & 7.04 & 606.30 \\
$\cotwo$-minimization  & 3.51& 4.89& 553.00\\
Weighted $w_f=2$ & 2.46& 5.24& 371.15\\
Flattening & 2.38& 5.60& 272.05 \\
\hline
\end{tabular}
\end{table}
This case has the highest power peaks ($P(s)=606.03$~kW) among the considered cases. Furthermore, we can clearly identify the periodicity in the schedule in Fig.~\ref{fig:bus:powerprofiles_noglobalcal_greedy}. Buses arrive starting in the afternoon and up until the early morning hours, creating a midnight-centric periodic pattern in the charging profile over the course of the week. It should be noted that the $\cotwo$ emission factor mirrors a rather sunny period, with high emission factors during the night, and low emission factors around noon. Combined with the pattern seen in the charging profile, most of the energy consumption in the uncontrolled scenario occurs during high-emission hours. Over the course of the week, this amounts to a total of 7.04~$\times10^{6}$~kg$\cotwo$eq in emissions.

The second scenario minimizes the associated $\cotwo$ emissions $C(s)$ using the minimum flow approach. Fig.~\ref{fig:bus:powerprofiles_noglobalcal_cotwo} clearly shows that the profile displays peaks during the low-emission periods throughout the day. However, the overall profile is also more segmented compared to the uncontrolled scenario as there are more blocks of high-power charging, also during more emission intensive periods. This results from the availability structure of the various charging jobs. Some amount of charging must take place during the night to complete the charging process before the buses depart. The total amount of $\cotwo$ associated with the schedule is 4.89$\times10^{6}$~kg$\cotwo$eq, which is a reduction of almost 31\% compared to the uncontrolled scenario. Furthermore, power peaks are reduced by around 9\% compared to the uncontrolled case.

The other extreme is depicted in Fig.~\ref{fig:bus:powerprofiles_noglobalcal_flat}. Here, the objective function is $F(s)$, flattening the power profile. The resulting profile can be partitioned into two aggregated charging modes of around 140~kW and 250~kW. The profile clearly displays periodicity with a frequency of about a day. Furthermore, the overal charging demand resulting from buses driving according to the weekend schedule is clearly lower than the demand resulting from buses driving on a week day. While, compared to the carbon-minimizing schedule, the associated emissions over the week increase to 5.60~$\times10^{6}$~kg$\cotwo$eq, which is lower than in the uncontrolled case. Moreover, the achieved objective value $F(s)$ for the scenario optimizing the flattening objective is only 54\% of the value of the uncontrolled case, and the power peak $P(s)$ reduces to 272.05~kW compared to 606.30~kW in the uncontrolled case, amounting to a reduction of 55\% in peak power compared to the uncontrolled case. Associated emissions of the schedule are about 20\% lower than in the uncontrolled case.

For the weighted objective function~\eqref{eq:weightedobjectiveDiscretized}, the results for $w_c = 1$ and $w_f = 2$ are illustratively presented in Fig.~\ref{fig:bus:powerprofiles_noglobalcal_weighted}. Even more so than when minimizing $C(s)$, we see the periods of high carbon intensity reflected in the power profile of the schedule. In particular, close to local minima and maxima of the emission factor plot, the power profile horizontally mirrors the curve of the emissions. The total emissions over the week are 5.24~$\times10^{6}$~kg$\cotwo$eq (26\% reduction compared to uncontrolled case), and the power peak is at 371.15~kW (39\% reduction compared to uncontrolled case). Those values are positioned between the extremes of purely steering for either carbon-minimization or flattening. 

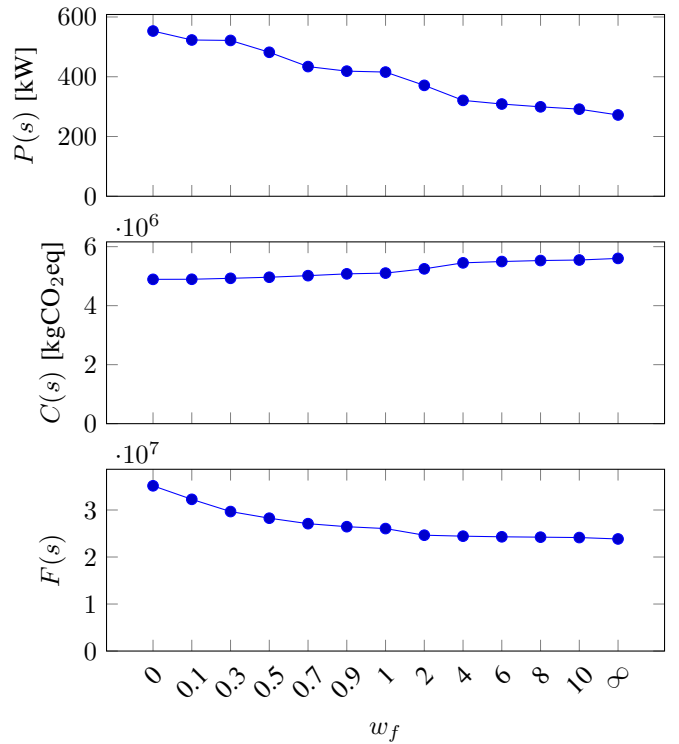
\begin{figure}
    \centering
    \hfill
        \begin{tikzpicture}
        \begin{axis}[
        height = 4cm,
        width = 9cm,
        ylabel={$P(s)$ [kW]},
        xticklabel = \empty,
        xtick={0,1,2,3,4,5,6,7,8,9,10,11,12}, 
        ymin=0,
        ]
        \addplot table [x expr=\coordindex, y=peak_size, col sep=comma] {cost_data.csv};
        \end{axis}
        \end{tikzpicture}
    \\
    \hfill
        \begin{tikzpicture}
        \begin{axis}[
            height = 4cm,
            width = 9cm,
            ylabel={$C(s)$ [kg$\cotwo {\text{eq}}$]},
            xticklabel = \empty,
            xtick={0,1,2,3,4,5,6,7,8,9,10,11,12},  
            ymin=0,
        ]
        \addplot table [x expr=\coordindex, y=cotwocost, col sep=comma] {cost_data.csv};
        \end{axis}
        \end{tikzpicture}
    \\
    \hfill
        \begin{tikzpicture}
        \begin{axis}[
            height = 4cm,
            width = 9cm,
            xlabel={$w_f$},
            ylabel={$F(s)$},
            xtick={0,1,2,3,4,5,6,7,8,9,10,11,12}, 
            xticklabels={0,0.1,0.3,0.5,0.7,0.9, 1, 2, 4, 6, 8, 10, $\infty$}, 
            xticklabel style={rotate=45},  
            ymin=0,
        ]
        \addplot table [x expr=\coordindex, y=flat_cost, col sep=comma] {cost_data.csv};
        \end{axis}
        \end{tikzpicture}
\caption{Metric results for the controlled scenarios for $w_c=1$ and various weights $w_f$.}
\label{fig:results_costs_metrics}
\end{figure}
To further illustrate the trade-off between the two objectives, Fig.~\ref{fig:results_costs_metrics} presents the three evaluation metrics for various weights $w_f$ and constant $w_c = 1$.
The extreme weights $w_f = 0$ and $w_f = \infty$ on the y-axis correspond to the cases that steer based on respectively objective functions $C(s)$ and $F(s)$.

\subsection{Bus Charging flexibility} \label{sec:results:baseloadflex}

In this section, we investigate how much flexibility is included in the bus charging case, as in the future it may be possible for a bus charging site to share a grid connection with other companies to further reduce stress on the grid and to free up capacity when it is scarce.

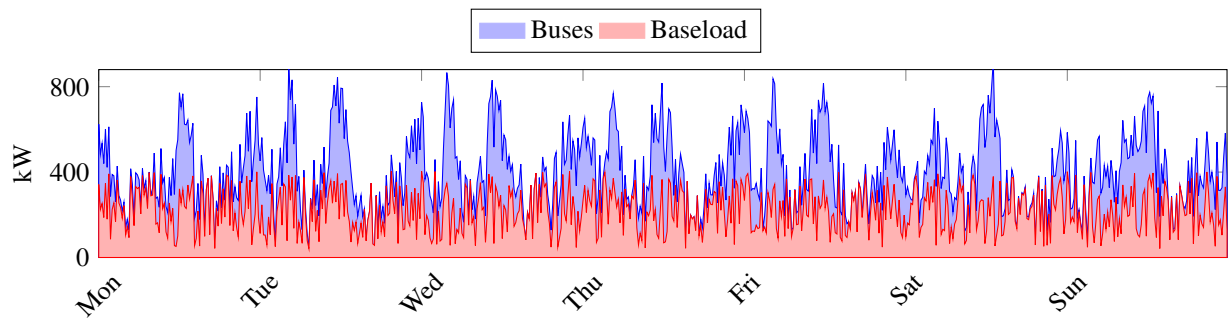
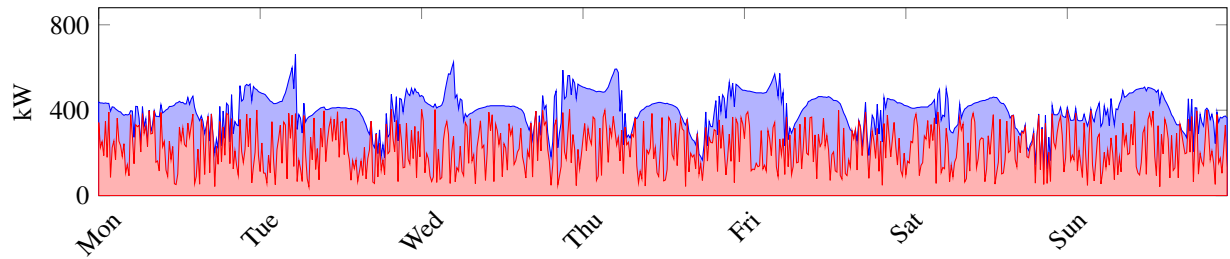
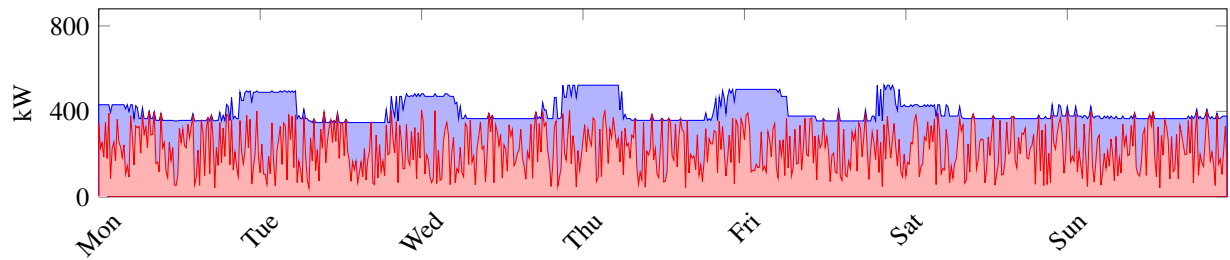
\begin{figure*}
    \begin{center}
        \ref{legend:bus:bigbaseload}
    \end{center}
    \centering
    \subfloat[The Power Profiles for CO$_2$ steering.]{\label{fig:power_profiles_rb2}
        \begin{tikzpicture}
            \pgfplotsset{
                scale only axis,
                xmin=0, xmax=671,
                height = 2.5cm, 
                width = 15cm,
                xtick={0,96,192,288,384,480,576},  
                xticklabels={Mon, Tue, Wed, Thu, Fri, Sat, Sun}, 
                xticklabel style={rotate=45}  
            }
            \begin{axis}[
                ytick = {0,100, 200},
                yticklabels = {0, 400, 800},
                ylabel = {kW},
                ymin = 0,
                ymax = 220
            ]
                \addplot[no marks, fill=blue!30, draw=blue] table[x expr=\coordindex, y index=1, col sep=comma, header=false]{power_profiles_random_baseload.csv} \closedcycle;
                \addplot[no marks, fill=red!30, draw=red] table[x expr=\coordindex, y index=0, col sep=comma]{combinedbaseload.csv} \closedcycle;
            \end{axis}
    \end{tikzpicture}
    }
    \\
    \subfloat[The Power Profiles for $w_c = 1$ and $w_f = 2$.]{\label{fig:power_profiles_rb3}
        \begin{tikzpicture}
            \pgfplotsset{
                scale only axis,
                xmin=0, xmax=671,
                height = 2.5cm, 
                width = 15cm,
                xtick={0,96,192,288,384,480,576}, 
                xticklabels={Mon, Tue, Wed, Thu, Fri, Sat, Sun},  
                xticklabel style={rotate=45}  
            }
            \begin{axis}[
                ytick = {0,100, 200},
                yticklabels = {0, 400, 800},
                ylabel = {kW},
                ymin = 0,
                ymax = 220
            ]
            \addplot[no marks, fill=blue!30, draw=blue] table[x expr=\coordindex, y index=7, col sep=comma, header=false]{power_profiles_random_baseload.csv} \closedcycle;
            \addplot[no marks, fill=red!30, draw=red] table[x expr=\coordindex, y index=0, col sep=comma]{combinedbaseload.csv} \closedcycle;
            \end{axis}
    \end{tikzpicture}
    }
    \\
    \subfloat[The Power Profiles for flattening steering.]{\label{fig:power_profiles_rb4} 
        \begin{tikzpicture}
            \pgfplotsset{
                scale only axis,
                xmin=0, xmax=671,
                height = 2.5cm, 
                width = 15cm,
                xtick={0,96,192,288,384,480,576},  
                xticklabels={Mon, Tue, Wed, Thu, Fri, Sat, Sun}, 
                xticklabel style={rotate=45}, 
                legend to name = legend:bus:bigbaseload,
                legend columns = 3
            }
            \begin{axis}[
                ytick = {0,100, 200},
                yticklabels = {0, 400, 800},
                ylabel = {kW},
                ymin = 0,
                ymax = 220
            ]
                \addplot[name path = A, color = blue, forget plot]table[x expr=\coordindex, y index=13, col sep=comma, header=false]{power_profiles_random_baseload.csv} \closedcycle; \label{plot_5}
                \addplot[name path = B, color = red, forget plot]table[x expr=\coordindex, y index=0, col sep=comma, header =false]{combinedbaseload.csv} 
                \closedcycle; \label{plot_two}
                \addplot[name path = C, white, forget plot] {0};

                \addplot[blue!30] fill between[of=A and B];
                \addlegendentry{Buses};

                \addplot[red!30] fill between[of = B and C];
                \addlegendentry{Baseload};
            \end{axis}
    \end{tikzpicture}
    }
    \caption{The power profiles with randomly sampled baseload.}
    \label{fig:power_profiles_rb}
\end{figure*}

For this, a dummy baseload is generated that consumes between $40$~kW and $400$~kW in each time interval according to a uniform probability distribution. This dummy baseload is then added to the real baseload recorded at the bus charging site. This profile is chosen to mimic unpredictable large fluctuations. By that, it is a good test case for flexibility as flattening this profile requires working around these large peaks. In Fig.~\ref{fig:power_profiles_rb} we present the results from the aforementioned steering methods for this case, i.e., the resulting profiles for the flattening scenario, the scenario with $\cotwo$ minimization, and the weighted scenario with $w_f=2$ are given. 

These figures show that the average power consumption has increased and that the fraction of the power that is allocated to the baseload also increased in comparison with Fig.~\ref{fig:bus:powerprofiles_noglobalcap}. That was to be expected as we added an additional baseload with an average of $220$~kW using a uniform distribution. 

In the case where we minimize the sum of squares $F(s)$ (see Fig.~\ref{fig:power_profiles_rb4}) the total profile is quite flat and displays similar periodic high and low power charging behaviour as in Fig.~\ref{fig:bus:powerprofiles_noglobalcal_flat}, despite the large magnitude and fluctuations of the added baseload.

To measure the flexibility included in the bus charging problem we take the largest peak size of the randomly generated baseload and the largest peak of the bus charging problem without this baseload added, where the schedule results from minimizing the sum of squares (so $w_f = \infty$). The sum of those peaks can then be compared to the largest peak of the case where both bus charging and baseload are considered during optimization. This indicates how much the worst case total peak can be reduced by coordinating both loads. 

For the use case of this paper, the baseload peak, the flattened bus charging peak, and the flattened peak when planning the bus charging around the baseload are respectively 404.2~kW, 272.0~kW and 522.4~kW.
When charging the buses independently of the fluctuating baseload, a worst-case peak of $404.2+272.0 = 676.2$~kW can occur. However, with coordination, this peak decreases to $522.4$~kW. The charging of the buses now only requires $118.2$~kW additional capacity, which corresponds to a decrease of $57\%$. We therefore conclude that it is worthwhile to coordinate bus charging sites with neighboring loads to decrease stress on grid assets.

\section{Conclusion} \label{bus:sec:conclusion}
In this work, we consider the scheduling of the charging activities at a charging site for electric buses from an aggregated perspective. The goal of this scheduling is to help the energy transition and the switch to renewable energy resources. To this end, we discussed three objective functions and a respective solver for each of them: $\cotwo$ minimization, sum of square minimization (flattening the aggregated power profile), and a weighted combination of the two. 
All presented methods were validated based on real-world data. The results show a clear trade-off between emission-intensity and flatness of a schedule for the weighted objective function. For example, the results indicate that purely optimizing for $\cotwo$ emission minimization reduces emissions by almost 31\%, but only reduces peak power by less than 9\%. 
Similarly, while purely optimizing for flatness resulted in a peak power reduction of 55\% compared to the uncontrolled case, it only reduced $\cotwo$ emissions by 20\%. In one of the weighted scenarios, power peaks and emissions were reduced by respectively 39\% and 26\%.

Furthermore, using a uniformly sampled random baseload, we showed the potential flexibility present in a bus charging site, indicating that it may be beneficial to integrate such a site into environments with other assets, e.g., data centres or large office buildings.

Future work may consider the bus-to-line matching in more detail, as in this work we only used a simple maximum matching model, minimizing the number of buses used. In this study, the effect of adding another bus or changing the considered bus-to-line matching on the objective function of the resulting schedule was left out of scope. Another direction may be to integrate the power capacity of the grid connection. While this is straightforward for the minimum cost flow method, its integration into the used \focs\ solver requires careful book-keeping, which as of today is incompatible with the available software implementation and therefore could not be validated in this study. 

\bibliographystyle{IEEEtran}
\bibliography{main}

\begin{thebibliography}{10}
\providecommand{\url}[1]{#1}
\csname url@samestyle\endcsname
\providecommand{\newblock}{\relax}
\providecommand{\bibinfo}[2]{#2}
\providecommand{\BIBentrySTDinterwordspacing}{\spaceskip=0pt\relax}
\providecommand{\BIBentryALTinterwordstretchfactor}{4}
\providecommand{\BIBentryALTinterwordspacing}{\spaceskip=\fontdimen2\font plus
\BIBentryALTinterwordstretchfactor\fontdimen3\font minus
  \fontdimen4\font\relax}
\providecommand{\BIBforeignlanguage}[2]{{%
\expandafter\ifx\csname l@#1\endcsname\relax
\typeout{** WARNING: IEEEtran.bst: No hyphenation pattern has been}%
\typeout{** loaded for the language `#1'. Using the pattern for}%
\typeout{** the default language instead.}%
\else
\language=\csname l@#1\endcsname
\fi
#2}}
\providecommand{\BIBdecl}{\relax}
\BIBdecl

\bibitem{hofstedeStroomnettenZijnOvervol2023}
{De Volkskrant}, ``{De stroomnetten zijn overvol, driekwart van Nederland heeft
  geen capaciteit voor nieuwe contracten voor grootschalige afname van
  elektriciteit},''
  https://www.volkskrant.nl/nieuws-achtergrond/de-stroomnetten-zijn-overvol-driekwart-van-nederland-heeft-geen-capaciteit-voor-nieuwe-contracten-voor-grootschalige-afname-van-elektriciteit{\textasciitilde}b52244e1/,
  Dec. 2023.

\bibitem{ToezichthouderACMWil2024}
{NOS}, ``{Toezichthouder ACM wil druk op stroomnet verminderen met pakket
  maatregelen},''
  https://nos.nl/artikel/2517302-toezichthouder-acm-wil-druk-op-stroomnet-verminderen-met-pakket-maatregelen,
  Apr. 2024.

\bibitem{perumalElectricBusPlanning2022}
S.~S.~G. Perumal, R.~M. Lusby, and J.~Larsen, ``Electric bus planning \&
  scheduling: {{A}} review of related problems and methodologies,''
  \emph{European Journal of Operational Research}, vol. 301, no.~2, pp.
  395--413, Sep. 2022.

\bibitem{2023KleinEVchargeSchedwithFlexibleServiceOperations}
\BIBentryALTinterwordspacing
P.~S. Klein and M.~Schiffer, ``Electric vehicle charge scheduling with flexible
  service operations,'' \emph{Transportation Science}, vol.~57, no.~6, pp.
  1605--1626, 2023. [Online]. Available:
  \url{https://doi.org/10.1287/trsc.2022.0272}
\BIBentrySTDinterwordspacing

\bibitem{zhengOnlineDistributedMPCBased2019}
Y.~Zheng, Y.~Song, D.~J. Hill, and K.~Meng, ``Online {{Distributed MPC-Based
  Optimal Scheduling}} for {{EV Charging Stations}} in {{Distribution
  Systems}},'' \emph{IEEE Transactions on Industrial Informatics}, vol.~15,
  no.~2, pp. 638--649, Feb. 2019.

\bibitem{shindeOptimalElectricVehicle2016}
P.~Shinde and K.~Swarup, ``Optimal {{Electric Vehicle}} charging schedule for
  demand side management,'' in \emph{2016 {{First International Conference}} on
  {{Sustainable Green Buildings}} and {{Communities}} ({{SGBC}})}, Dec. 2016,
  pp. 1--6.

\bibitem{alonsoOptimalChargingScheduling2014}
M.~Alonso, H.~Amaris, J.~Germain, and J.~Galan, ``Optimal charging scheduling
  of electric vehicles in smart grids by heuristic algorithms,''
  \emph{Energies}, vol.~7, no.~4, pp. 2449--2475, 2014.

\bibitem{adetunjiTwoTailedPricingScheme2024}
K.~E. Adetunji, I.~W. Hofsajer, A.~M. {Abu-Mahfouz}, and L.~Cheng, ``A
  {{Two-Tailed Pricing Scheme}} for {{Optimal EV Charging Scheduling Using
  Multiobjective Reinforcement Learning}},'' \emph{IEEE Transactions on
  Industrial Informatics}, vol.~20, no.~3, pp. 3361--3370, Mar. 2024.

\bibitem{zhuang_stochastic_2021}
\BIBentryALTinterwordspacing
P.~Zhuang and H.~Liang, ``Stochastic {Energy} {Management} of {Electric} {Bus}
  {Charging} {Stations} {With} {Renewable} {Energy} {Integration} and {B2G}
  {Capabilities},'' \emph{IEEE Transactions on Sustainable Energy}, vol.~12,
  no.~2, pp. 1206--1216, Apr. 2021, conference Name: IEEE Transactions on
  Sustainable Energy. [Online]. Available:
  \url{https://ieeexplore.ieee.org/document/9266099/?arnumber=9266099}
\BIBentrySTDinterwordspacing

\bibitem{liCooperativeOptimizationBus2023}
P.~Li, M.~Jiang, Y.~Zhang, and Y.~Zhang, ``Cooperative {{Optimization}} of
  {{Bus Service}} and {{Charging Schedules}} for a {{Fast-Charging Battery
  Electric Bus Network}},'' \emph{IEEE Transactions on Intelligent
  Transportation Systems}, vol.~24, no.~5, pp. 5362--5375, May 2023.

\bibitem{panOptimizationElectricBus2020}
J.~Pan, X.~Wu, Q.~Feng, and Y.~Ji, ``Optimization of {{Electric Bus Charging
  Station Considering Energy Storage System}},'' in \emph{2020 8th
  {{International Conference}} on {{Power Electronics Systems}} and
  {{Applications}} ({{PESA}})}.\hskip 1em plus 0.5em minus 0.4em\relax Hong
  Kong, China: IEEE, Dec. 2020.

\bibitem{basmaOptimizationBatteryElectric2020}
H.~Basma, C.~Mansour, M.~Nemer, M.~Haddad, and P.~Stabat, ``Optimization of
  {{Battery Electric Bus Charging}} under {{Varying Operating Conditions}},''
  in \emph{2020 {{IEEE Vehicle Power}} and {{Propulsion Conference}}
  ({{VPPC}})}.\hskip 1em plus 0.5em minus 0.4em\relax Gijon, Spain: IEEE, Nov.
  2020, pp. 1--6.

\bibitem{chengweizhangChargingScheduleOptimization2019}
{Chengwei Zhang}, ``Charging {{Schedule Optimization}} of {{Electric Bus
  Charging Station}} considering {{Departure Timetable}},'' in \emph{8th
  {{Renewable Power Generation Conference}} ({{RPG}} 2019)}.\hskip 1em plus
  0.5em minus 0.4em\relax Shanghai, China: {Institution of Engineering and
  Technology}, 2019.

\bibitem{leouOptimalChargingSchedule2017}
R.-C. Leou and J.-J. Hung, ``Optimal {{Charging Schedule Planning}} and
  {{Economic Analysis}} for {{Electric Bus Charging Stations}},''
  \emph{Energies}, vol.~10, no.~4, p. 483, Apr. 2017.

\bibitem{zhouCollaborativeOptimizationVehicle2020}
G.-J. Zhou, D.-F. Xie, X.-M. Zhao, and C.~Lu, ``Collaborative {{Optimization}}
  of {{Vehicle}} and {{Charging Scheduling}} for a {{Bus Fleet Mixed With
  Electric}} and {{Traditional Buses}},'' \emph{IEEE Access}, vol.~8, pp.
  8056--8072, 2020.

\bibitem{abdelwahedEvaluatingOptimizingOpportunity2020}
A.~Abdelwahed, P.~L. {van den Berg}, T.~Brandt, J.~Collins, and W.~Ketter,
  ``Evaluating and {{Optimizing Opportunity Fast-Charging Schedules}} in
  {{Transit Battery Electric Bus Networks}},'' \emph{Transportation Science},
  vol.~54, no.~6, pp. 1601--1615, Nov. 2020.

\bibitem{heOptimalChargingScheduling2020a}
Y.~He, Z.~Liu, and Z.~Song, ``Optimal charging scheduling and management for a
  fast-charging battery electric bus system,'' \emph{Transportation Research
  Part E: Logistics and Transportation Review}, vol. 142, p. 102056, Oct. 2020.

\bibitem{bagherinezhadSpatioTemporalElectricBus2020}
A.~Bagherinezhad, A.~D. Palomino, B.~Li, and M.~Parvania, ``Spatio-{{Temporal
  Electric Bus Charging Optimization With Transit Network Constraints}},''
  \emph{IEEE Transactions on Industry Applications}, vol.~56, no.~5, pp.
  5741--5749, Sep. 2020.

\bibitem{jahicChargingScheduleLoad2019}
A.~Jahic, M.~Eskander, and D.~Schulz, ``Charging {{Schedule}} for {{Load Peak
  Minimization}} on {{Large-Scale Electric Bus Depots}},'' \emph{Applied
  Sciences}, vol.~9, no.~9, p. 1748, Jan. 2019.

\bibitem{2023WinschermannFOCSArxiv}
\BIBentryALTinterwordspacing
L.~Winschermann, M.~E.~T. Gerards, A.~Antoniadis, G.~Hoogsteen, and J.~Hurink,
  ``Relating electric vehicle charging to speed scaling with job-specific speed
  limits,'' 2023, currently under peer-review. [Online]. Available:
  \url{https://arxiv.org/abs/2309.06174}
\BIBentrySTDinterwordspacing

\bibitem{2024nednlData}
\BIBentryALTinterwordspacing
{Nationaal Energie Dashboard}. (2024) Totale elektriciteitsproductie. Data in
  used granularity received via e-mail. [Online]. Available:
  \url{https://ned.nl/nl/dataportaal/energie-productie/elektriciteit/totale-elektriciteitsproductie}
\BIBentrySTDinterwordspacing

\end{thebibliography}

\end{document}